\documentclass[12pt]{article}
\usepackage{amssymb}

\newtheorem{thm}{Theorem}[section]

\newtheorem{proposition}[thm]{Proposition}
\newtheorem{definition}[thm]{Definition}
\newtheorem{remark}[thm]{Remark}

\begin{document}

\title{On the algebra ${\mathcal{A}}_{\hbar,\eta}(osp(2|2)^{(2)})$ and free boson
representations}

\author{Niall MacKay${}^{1}$\thanks{e-mail: nm15@york.ac.uk}
\hspace{0.5cm} and \hspace{0.5cm} Liu Zhao${}^{1,2}$\thanks{e-mail: lzhao@phy.nwu.edu.cn}
{}\thanks{%
Royal Society visiting fellow}\\
${}^1$ Department of Mathematics, University of York, York YO10 5DD, UK\\
${}^2$ Institute of Modern Physics, Northwest University, Xian 710069, China%
}
\date{}
\maketitle

\begin{abstract}
A two-parameter quantum deformation of the affine Lie super
algebra $osp(2|2)^{(2)}$ is introduced and studied in some detail.
This algebra is the first example associated with nonsimply-laced
and twisted root systems of a quantum current algebra with the
structure of a so-called infinite Hopf family of (super)algebras.
A representation of this algebra at $c=1$ is realized in the
product Fock space of two commuting sets of Heisenberg algebras.
\end{abstract}

\section{\protect\vspace{1pt}Introduction}

The study of two-parameter deformations of affine Lie
(super)algebras has turned out to be quite fruitful in the last
few years. These algebras are in a sense deformations of the
standard quantum affine algebras or Yangian doubles which all have
the structure of quasi-triangular Hopf algebras. In contrast, the
two-parameter deformations are not Hopf algebras, but rather have
twisted Hopf structures, one of which is the Drinfeld twist of
Hopf algebras or quasi-Hopf algebras, while the other is the
infinite Hopf family of (super)algebras. So far, the relationship
between the two generalized co-structures remains ill-understood.

\vspace{1pt}

{}From the physical point of view, the above two classes of
two-parameter deformations appear in different contexts.
Quasi-Hopf algebras \cite{D4} occur in the study of symmetries of
face- and vertex-type models of statistical mechanics
\cite{jimbot}, and are closely related to the face-type Boltzmann
weights (Yang-Baxter R-matrix), while infinite Hopf families of
(super)algebras \cite{f,hou2} occur only in the representation
theory of quantum deformed Virasoro and W-algebras \cite{FF} (and
act as the algebra of screening currents), which in turn are
algebras characterizing the dynamical symmetries of certain
massive integrable quantum field theories \cite{KLP}.

\vspace{1pt}

Despite the great deal of work that has been done on these
two-parameter deformations, many problems remain unsolved. In
particular, for the second class of two-parameter deformations
(the infinite Hopf families) nothing has been said concerning root
systems of nonsimply-laced and/or twisted types. For
super root systems, the only case which has been considered is the case of $%
osp(1|2)^{(1)}$ \cite{ZD}. In this article, therefore, we aim to
provide more concrete examples of this kind, and in particular the
trigonometric two-parameter deformation of $osp(2|2)^{(2)}$. This
algebra is based on a root system which is simultaneously
non-simply laced, twisted and super.

\vspace{1pt}

\section{Definition and structure of the algebra $\mathcal{A}_{\hbar ,\eta
}(osp(2|2)^{(2)})$}

\subsection{\protect\vspace{1pt}Definition of $\mathcal{A}_{\hbar ,\eta
}(osp(2|2)^{(2)}).$}

We start with the definition of the algebra $\mathcal{A}_{\hbar
,\eta }(osp(2|2)^{(2)})$. The notation follows that of $%
\mathcal{A}_{\hbar ,\eta }(\hat{g})$ \cite{KLP,f,hou2} and $\mathcal{A}_{\hbar ,\eta
}(osp(1|2)^{(1)})$ \cite{ZD}.

\begin{definition} \label{prop1}
The algebra $\mathcal{A}_{\hbar ,\eta }(osp(2|2)^{(2)})$, considered as a
continuously distributed current superalgebra, is a $Z_{2}$ graded
associative algebra over $C$ generated by the currents $E(u)$, $F(u)$, $%
H^{\pm }(u)$, the central element $c$ and unit $1$ with parities $\pi
\lbrack E(u)\rbrack =\pi \lbrack F(u)\rbrack =1,$ $\pi \lbrack H^{\pm
}(u)\rbrack =\pi \lbrack c\rbrack =\pi \lbrack 1\rbrack =0$ and generating
relations
\begin{eqnarray}
E(u)E(v) &=&-\frac{\cosh \lbrack \pi \eta (u-v+i\hbar )\rbrack }{\cosh
\lbrack \pi \eta (u-v-i\hbar )\rbrack }E(v)E(u), \\
F(u)F(v) &=&-\frac{\cosh \lbrack \pi \eta ^{\prime }(u-v-i\hbar )\rbrack }{%
\cosh \lbrack \pi \eta ^{\prime }(u-v+i\hbar )\rbrack }F(v)F(u), \\
H^{\pm }(u)E(v) &=&\frac{\cosh \lbrack \pi \eta (u-v+i\hbar \pm i\hbar
c/4)\rbrack }{\cosh \lbrack \pi \eta (u-v-i\hbar \pm i\hbar c/4)\rbrack }%
E(v)H^{\pm }(u), \\
H^{\pm }(u)F(v) &=&\frac{\cosh \lbrack \pi \eta ^{\prime }(u-v-i\hbar \mp
i\hbar c/4)\rbrack }{\cosh \lbrack \pi \eta ^{\prime }(u-v+i\hbar \mp i\hbar
c/4)\rbrack }F(v)H^{\pm }(u), \\
H^{\pm }(u)H^{\pm }(v) &=&\frac{\cosh \lbrack \pi \eta (u-v+i\hbar )\rbrack
}{\cosh \lbrack \pi \eta (u-v-i\hbar )\rbrack }\frac{\cosh \lbrack \pi \eta
^{\prime }(u-v-i\hbar )\rbrack }{\cosh \lbrack \pi \eta ^{\prime
}(u-v+i\hbar )\rbrack }H^{\pm }(v)H^{\pm }(u), \\
H^{+}(u)H^{-}(v) &=&\frac{\cosh \lbrack \pi \eta (u-v+i\hbar +i\hbar
c/2)\rbrack }{\cosh \lbrack \pi \eta (u-v-i\hbar +i\hbar c/2)\rbrack }\frac{%
\cosh \lbrack \pi \eta ^{\prime }(u-v-i\hbar -i\hbar c/2)\rbrack }{\cosh
\lbrack \pi \eta ^{\prime }(u-v+i\hbar -i\hbar c/2)\rbrack }H^{-}(v)H^{+}(u),
\\
\{E(u),F(v)\} &=&\frac{2\pi }{\hbar }\left[ \delta (u-v-i\hbar
c/2)H^{+}(u-i\hbar c/4)-\delta (u-v+i\hbar c/2)H^{-}(v-i\hbar c/4)\right] ,
\label{ant1}
\end{eqnarray}
where
\[
\frac{1}{\eta ^{\prime }}-\frac{1}{\eta }=\hbar c,
\]
and $\hbar $ and $\eta $ are generic deformation parameters.\hfill $\square $
\end{definition}

For later reference, we denote the subalgebras generated respectively by the
currents $\{E(u)\}$ and $\{F(u)\}$ as $\mathcal{N}_{\pm }\lbrack \mathcal{A}%
_{\hbar ,\eta }(osp(2|2)^{(2)})\rbrack $.

It is interesting to compare the generating relations of the algebra $%
\mathcal{A}_{\hbar ,\eta }(osp(2|2)^{(2)})$ and those of the quantum affine
super algebra $U_{q}(osp(2|2)^{(2)})$. The latter has the generating
relations \cite{YZ}
\begin{eqnarray}
X^{+}(z)X^{+}(w) &=&-\frac{zq+w}{z+wq}X^{+}(w)X^{+}(z), \label{uq}\\
X^{-}(z)X^{-}(w) &=&-\frac{z+wq}{zq+w}X^{-}(w)X^{-}(w), \\
\psi ^{\pm }(z)X^{+}(w) &=&\frac{z_{\pm }q+w}{z_{\pm }+wq}X^{+}(w)\psi ^{\pm
}(z), \\
\psi ^{\pm }(z)X^{-}(w) &=&\frac{z_{\mp }+wq}{z_{\mp }q+w}X^{-}(w)\psi ^{\pm
}(z), \\
\psi ^{\pm }(z)\psi ^{\pm }(w) &=&\psi ^{\pm }(w)\psi ^{\pm }(z), \\
\psi ^{+}(z)\psi ^{-}(w) &=&\frac{z_{+}q+w_{-}}{z_{+}+w_{-}q}\frac{%
z_{-}+w_{+}q}{z_{-}q+w_{+}}\psi ^{-}(w)\psi ^{+}(z), \\
\{X^{+}(z),X^{-}(w)\} &=&\frac{1}{(q-q^{-1})zw}\left[ \delta (\frac{w}{z}%
\gamma )\psi ^{+}(z_{-})-\delta (\frac{w}{z}\gamma ^{-1})\psi
^{-}(w_{-})\right] , \label{ant2}
\end{eqnarray}
where\footnote{%
In the original presentation of $U_{q}(osp(2|2)^{(2)})$ in
\cite{YZ}, the element $\gamma $ was written as $q^{c}$. However,
to avoid confusion with the central element $c$ of the algebra
$\mathcal{A}_{\hbar ,\eta }(osp(2|2)^{(2)})$, we intentionally
rename it $\gamma $, as is usual in ordinary quantum affine
algebras.} $z_{\pm }=z\gamma ^{\pm 1/2}$ . Notice that the $\delta
$-functions appearing in (\ref{ant1}) and (\ref{ant2}) are
supported differently, the former at $0$ ({\em i.e.\ }the standard
Dirac $\delta$-function), the latter at $1$ (\,$\delta(z) \equiv
\sum_{n \in Z} z^n$ as a formal power series). We also denote the
subalgebras of $U_{q}(osp(2|2)^{(2)})$ generated respectively by
$X^{+}(z)$ and $X^{-}(z)$ by $\mathcal{N}_{\pm }\lbrack
U_{q}(osp(2|2)^{(2)})\rbrack $.

The following two propositions justify the similarities between
the two
algebras $\mathcal{A}_{\hbar ,\eta }(osp(2|2)^{(2)})$ and $%
U_{q}(osp(2|2)^{(2)})$: first,

\begin{proposition} \label{prop2}
There are algebra homomorphisms $\rho ^{+}:\mathcal{N}_{+}\lbrack \mathcal{A}%
_{\hbar ,\eta }(osp(2|2)^{(2)})\rbrack \rightarrow \mathcal{N}_{+}\lbrack
U_{q}(osp(2|2)^{(2)})\rbrack $, $\rho ^{-}:\mathcal{N}_{-}\lbrack \mathcal{A}%
_{\hbar ,\eta }(osp(2|2)^{(2)})\rbrack \rightarrow
\mathcal{N}_{-}\lbrack U_{q^{\prime }}(osp(2|2)^{(2)})\rbrack $,
where under $\rho ^{\pm }$, the parameters behave as
\[
\begin{array}{l}
\rho ^{+}(e^{2\pi \eta u})=z, \\ \rho ^{+}(e^{2\pi i\eta \hbar
})=q,
\end{array}
\]
and
\[
\begin{array}{l}
\rho ^{-}(e^{2\pi \eta ^{\prime }u})=z, \\ \rho ^{-}(e^{2\pi i\eta
^{\prime }\hbar })=q^{\prime },
\end{array}
\]
respectively.\hfill $\square $
\end{proposition}

Recalling that $\eta $ and $\eta ^{\prime }$ are different only
when $c\neq 0$, we also have

\begin{proposition} \label{prop3}
There is an algebra homomorphism between $\mathcal{A}_{\hbar ,\eta
}(osp(2|2)^{(2)})$ at $c=0$ and $U_{q}(osp(2|2)^{(2)})$ at $\gamma
=1$:
\begin{eqnarray*}
\mathcal{E}:\mathcal{A}_{\hbar ,\eta }(osp(2|2)^{(2)}) &\rightarrow
&U_{q}(osp(2|2)^{(2)}), \\
E(u) &\longmapsto &\sqrt{2}zX^{+}(z), \\
F(u) &\longmapsto &\sqrt{2}zX^{-}(z), \\
\frac{2\pi }{\hbar }H^{\pm }(u) &\longmapsto &\frac{1}{q-q^{-1}}\psi ^{\pm
}(z),
\end{eqnarray*}
where $z=e^{2\pi \eta \hbar u},$ $q=e^{2\pi i\eta \hbar }$.\hfill $\square $
\end{proposition}

Proposition \ref{prop2} indicates that the algebra
$\mathcal{A}_{\hbar ,\eta }(osp(2|2)^{(2)})$ is actually an
interpolation between (Borel subalgebras of) two standard quantum
affine algebras $U_{q}(osp(2|2)^{(2)})$ and $U_{q^{\prime
}}(osp(2|2)^{(2)})$ with different deformation parameters, while
Proposition \ref{prop3} further states that at $c=0$, the algebra
$\mathcal{A}_{\hbar ,\eta }(osp(2|2)^{(2)})$ degenerates into
$U_{q}(osp(2|2)^{(2)})$ at $\gamma =1$.

\subsection{Co-structure}

As expected, this algebra possesses the structure of an infinite Hopf family
of super algebras, whose definition can be found in \cite{ZD} (see
also \cite{f, hou2}). In fact, if we denote $\mathcal{A}_{n}=\mathcal{A}
_{\hbar ,\eta ^{(n)}}(osp(2|2)^{(2)})_{c_{n}}$ where $\eta ^{(n)}$ is
defined iteratively via $\frac{1}{\eta ^{(n+1)}}-\frac{1}{\eta ^{(n)}}=
\hbar c_{n}$ starting from $\eta ^{(1)}=\eta $ and taking $%
c_{n}\in Z\backslash Z_{-}$, we can define the following co-structures over
the family of algebras $\{\mathcal{A}_{n}$, $n\in Z\}$:

\begin{itemize}
\item  the comultiplications $\Delta _{n}^{\pm }$ (algebra homomorphisms $%
\Delta _{n}^{+}:\mathcal{A}_{n}\rightarrow \mathcal{A}_{n}\otimes \mathcal{A}%
_{n+1}$\emph{, }$\Delta _{n}^{-}:\mathcal{A}_{n}\rightarrow \mathcal{A}%
_{n-1}\otimes \mathcal{A}_{n}$):
\begin{eqnarray*}
\Delta _{n}^{+}c_{n}\!\!\!\! &=&\!\!\!\!c_{n}+c_{n+1}, \\
\Delta _{n}^{-}c_{n}\!\!\!\! &=&\!\!\!\!c_{n-1}+c_{n}, \\
\Delta _{n}^{+}H^{+}(u;\eta ^{(n)})\!\!\!\! &=&\!\!\!\!H^{+}(u+\frac{i\hbar
c_{n+1}}{4};\eta ^{(n)})\otimes H^{+}(u-\frac{i\hbar c_{n}}{4};\eta
^{(n+1)}), \\
\Delta _{n}^{-}H^{+}(u;\eta ^{(n)})\!\!\!\! &=&\!\!\!\!H^{+}(u+\frac{i\hbar
c_{n}}{4};\eta ^{(n-1)})\otimes H^{+}(u-\frac{i\hbar c_{n-1}}{4};\eta
^{(n)}), \\
\Delta _{n}^{+}H^{-}(u;\eta ^{(n)})\!\!\!\! &=&\!\!\!\!H^{-}(u-\frac{i\hbar
c_{n+1}}{4};\eta ^{(n)})\otimes H^{-}(u+\frac{i\hbar c_{n}}{4};\eta
^{(n+1)}), \\
\Delta _{n}^{-}H^{-}(u;\eta ^{(n)})\!\!\!\! &=&\!\!\!\!H^{-}(u-\frac{i\hbar
c_{n}}{4};\eta ^{(n-1)})\otimes H^{-}(u+\frac{i\hbar c_{n-1}}{4};\eta
^{(n)}), \\
\Delta _{n}^{+}E(u;\eta ^{(n)})\!\!\!\! &=&\!\!\!\!E(u;\eta ^{(n)})\otimes
1+H^{-}(u+\frac{i\hbar c_{n}}{4};\eta ^{(n)})\otimes E(u+\frac{i\hbar c_{n}}{%
2};\eta ^{(n+1)}), \\
\Delta _{n}^{-}E(u;\eta ^{(n)})\!\!\!\! &=&\!\!\!\!E(u;\eta ^{(n-1)})\otimes
1+H^{-}(u+\frac{i\hbar c_{n-1}}{4};\eta ^{(n-1)})\otimes E(u+\frac{i\hbar
c_{n-1}}{2};\eta ^{(n)}), \\
\Delta _{n}^{+}F(u;\eta ^{(n)})\!\!\!\! &=&\!\!\!\!1\otimes F(u;\eta
^{(n+1)})+F(u+\frac{i\hbar c_{n+1}}{2};\eta ^{(n)})\otimes H^{+}(u+\frac{%
i\hbar c_{n+1}}{4};\eta ^{(n+1)}), \\
\Delta _{n}^{-}F(u;\eta ^{(n)})\!\!\!\! &=&\!\!\!\!1\otimes F(u;\eta
^{(n)})+F(u+\frac{i\hbar c_{n}}{2};\eta ^{(n-1)})\otimes H^{+}(u+\frac{%
i\hbar c_{n}}{4};\eta ^{(n)});
\end{eqnarray*}

\item  the counits $\epsilon _{n}$ (algebra homomorphism $\epsilon _{n}:%
\mathcal{A}_{n}\rightarrow C):$
\begin{eqnarray*}
\epsilon _{n}(c_{n})\!\!\!\! &=&\!\!\!\!0, \\
\epsilon _{n}(1_{n})\!\!\!\! &=&\!\!\!\!1, \\
\epsilon _{n}(H_{i}^{\pm }(u;\eta ^{(n)}))\!\!\!\! &=&\!\!\!\!1, \\
\epsilon _{n}(E_{i}(u;\eta ^{(n)}))\!\!\!\! &=&\!\!\!\!0, \\
\epsilon _{n}(F_{i}(u;\eta ^{(n)}))\!\!\!\! &=&\!\!\!\!0;
\end{eqnarray*}

\item  the antipodes $S_{n}^{\pm }$ (algebra anti-homomorphisms $S_{n}^{\pm
}:\mathcal{A}_{n}\rightarrow \mathcal{A}_{n\pm 1}$):

\begin{eqnarray*}
S_{n}^{\pm }(c_{n})\!\!\!\! &=&\!\!\!\!-c_{n\pm 1}, \\
S_{n}^{\pm }(H^{\pm }(u;\eta ^{(n)}))\!\!\!\! &=&\!\!\!\!\lbrack H^{\pm
}(u;\eta ^{(n\pm 1)})\rbrack ^{-1}, \\
S_{n}^{\pm }(E(u;\eta ^{(n)}))\!\!\!\! &=&\!\!\!\!-H^{-}(u-\frac{i\hbar
c_{n\pm 1}}{4};\eta ^{(n\pm 1)})^{-1}E(u-\frac{i\hbar c_{n\pm 1}}{2};\eta
^{(n\pm 1)}), \\
S_{n}^{\pm }(F(u;\eta ^{(n)}))\!\!\!\! &=&\!\!\!\!-F(u-\frac{i\hbar c_{n\pm
1}}{2};\eta ^{(n\pm 1)})H^{+}(u-\frac{i\hbar c_{n\pm 1}}{4};\eta ^{(n\pm
1)})^{-1},
\end{eqnarray*}
\end{itemize}

\noindent where $\otimes $ stands for the direct super product defined by
\[
(A\otimes B)(C\otimes D)=(-1)^{\pi \lbrack B\rbrack \pi \lbrack C\rbrack
}AB\otimes CD
\]
for homogeneous elements $A,B,C,D$. It is a trivial (but tedious)
exercise to check that these structures satisfy all the defining
axioms of an infinite Hopf family of super algebras,

\begin{itemize}
\item  $(\epsilon _{n}\otimes id_{n+1})\circ \Delta _{n}^{+}=\tau
_{n}^{+},~(id_{n-1}\otimes \epsilon _{n})\circ \Delta _{n}^{-}=\tau _{n}^{-}$%
\emph{\ \hfill }

\item  $m_{n+1}\circ (S_{n}^{+}\otimes id_{n+1})\circ \Delta
_{n}^{+}=\epsilon _{n+1}\circ \tau _{n}^{+},~m_{n-1}\circ (id_{n-1}\otimes
S_{n}^{-})\circ \Delta _{n}^{-}=\epsilon _{n-1}\circ \tau _{n}^{-}$\ \hfill

\item  $(\Delta _{n}^{-}\otimes id_{n+1})\circ \Delta
_{n}^{+}=(id_{n-1}\otimes \Delta _{n}^{+})\circ \Delta _{n}^{-}$\ \hfill
\end{itemize}
\noindent where $m_{n}$\ is the (super)multiplication for $\mathcal{A}_{n}$ and $\tau
_{n}^{\pm }$ are algebra shift morphisms $\tau _{n}^{\pm }:\mathcal{A}%
_{n}\rightarrow \mathcal{A}_{n\pm 1}$which obeys $\tau _{n+1}^{-}\tau
_{n}^{+}=id_{n}=$ $\tau _{n-1}^{+}\tau _{n}^{-}$.

\vspace{1pt}

The operations $\Delta _{n}^{\pm }$ are related to each other by
the shift morphisms:
\begin{eqnarray*}
\Delta _{n}^{-} &=&(\tau _{n}^{-}\otimes \tau _{n+1}^{-})\circ \Delta
_{n}^{+}, \\
\Delta _{n}^{+} &=&(\tau _{n-1}^{+}\otimes \tau _{n}^{+})\circ \Delta
_{n}^{-}.
\end{eqnarray*}
Thus the easily-observed co-commutativity between the two
co-multiplications
\[
(\Delta _{n}^{-}\otimes id_{n+1})\circ \Delta _{n}^{+}=(id_{n-1}\otimes
\Delta _{n}^{+})\circ \Delta _{n}^{-}
\]
can be rewritten in terms of only one of the two
co-multiplications, and turns out to become a statement of the
non-coassociativity of the co-multiplications:
\[
\lbrack ((\tau _{n}^{-}\otimes \tau _{n+1}^{-})\circ \Delta _{n}^{+})\otimes
id_{n+1}\rbrack \circ \Delta _{n}^{+}=(id_{n-1}\otimes \Delta _{n}^{+})\circ
((\tau _{n}^{-}\otimes \tau _{n+1}^{-})\circ \Delta _{n}^{+}),
\]
\[
(\Delta _{n}^{-}\otimes id_{n+1})\circ ((\tau _{n-1}^{+}\otimes \tau
_{n}^{+})\circ \Delta _{n}^{-})=\lbrack id_{n-1}\otimes ((\tau
_{n-1}^{+}\otimes \tau _{n}^{+})\circ \Delta _{n}^{-})\rbrack \circ \Delta
_{n}^{-}.
\]
Notice that these twisted co-associativity conditions are
different from that of the Drinfeld twists. However, the effects
of these two different kinds of twists are the same: they all
allow one to construct fused (tensor product) representations for
the algebras under investigation, although the co-structures are
not co-associative.

\vspace{1pt}

Now recall that definition \ref{prop1} defines the algebra $%
\mathcal{A}_{\hbar ,\eta }(osp(2|2)^{(2)})$ only as a formal
algebra, in the sense that all currents thus defined are actually
only distributions. To assign precise meaning to the algebra
$\mathcal{A}_{\hbar ,\eta }(osp(2|2)^{(2)})$ we need to specify
the actual generators and relations, and this can be done only
separately for two distinct cases $c=0$ and $c\neq
0$ (as in the case of $\mathcal{A}_{\hbar ,\eta }(\widehat{sl_{2}})$ \cite{KLP} and $%
\mathcal{A}_{\hbar ,\eta }(\hat{g})$ \cite{f}). For details, the
reader is directed to Khoroshkin {\em et al.\ }\cite{KLP} in the
$\widehat{sl}_2$ case. The present case is in complete analogy.

\section{Representation theory}

\subsection{Case $c=0$}

Recall that for $c=0,$ there is an algebra homomorphism between the algebras
$\mathcal{A}_{\hbar ,\eta }(osp(2|2)^{(2)})$ and $U_{q}(osp(2|2)^{(2)})$ for
$q=e^{2\pi i\eta \hbar }$. Thus the evaluation representation of $%
U_{q}(osp(2|2)^{(2)})$ presented in \cite{YZ} can be extended into
an evaluation representation of $\mathcal{A}_{\hbar ,\eta }(osp(2|2)^{(2)})$
in terms of the evaluation homomorphism $\mathcal{E}$. This evaluation
representation justifies the relationship between the algebra $\mathcal{A}%
_{\hbar ,\eta }(osp(2|2)^{(2)})$ and the root system of type $osp(2|2)^{(2)}$%
.

\vspace{1pt}
\subsection{Case $c=1$ and structure of the Fock space}

As usual, the tool we need to construct a representation of $\mathcal{A}%
_{\hbar ,\eta }(osp(2|2)^{(2)})$ at $c=1$ is the free boson
realization. Throughout this subsection we have $1/\eta ^{\prime
}=1/\eta +\hbar $.

\vspace{1pt} Define the Heisenberg algebras $\mathcal{H}_{\alpha
}$, $\mathcal{H}_{\beta } $ respectively by
\begin{eqnarray*}
\lbrack \alpha (\lambda ),\alpha (\mu )\rbrack &=&A(\lambda )\delta (\lambda
+\mu ), \\
\lbrack \beta (\lambda ),\beta (\mu )\rbrack &=&B(\lambda )\delta (\lambda
+\mu ), \\
\lbrack \alpha (\lambda ),\beta (\mu )\rbrack &=&0, \\
(\lambda &\neq &\mu )
\end{eqnarray*}
where $A(\lambda )$ and $B(\lambda )$ are given as
\begin{eqnarray*}
A(\lambda ) &=&\frac{\lambda }{4\cosh \frac{\hbar \lambda }{2}+\left(
\mathrm{csch}\frac{\lambda }{2\eta }-\mathrm{csch}\frac{\lambda }{2\eta
^{\prime }}\right) \sinh \hbar \lambda +2}, \\
B(\lambda ) &=&\frac{\lambda \left( \left( 1+\mathrm{csch}\frac{\lambda }{%
2\eta }\sinh \hbar \lambda \right) \left( 1-\mathrm{csch}\frac{\lambda }{%
2\eta ^{\prime }}\sinh \hbar \lambda \right) -4\cosh ^{2}\frac{\hbar \lambda
}{2}\right) }{4\cosh \frac{\hbar \lambda }{2}+\left( \mathrm{csch}\frac{%
\lambda }{2\eta }-\mathrm{csch}\frac{\lambda }{2\eta ^{\prime }}\right)
\sinh \hbar \lambda +2},
\end{eqnarray*}
both of which are antisymmetric as $\lambda \rightarrow -\lambda $
and regular as $\lambda \rightarrow 0$. In fact we can easily
check that
\begin{eqnarray*}
A(\lambda ) &=&-A(-\lambda ), \\
B(\lambda ) &=&-B(-\lambda ),
\end{eqnarray*}
and
\begin{eqnarray*}
A(\lambda ) &\sim &\frac{\lambda }{2(\eta -\eta ^{\prime })\hbar +6}%
+O(\lambda ^{3}), \\
B(\lambda ) &\sim &\frac{\lbrack (1+2\hbar \eta )(1-2\hbar \eta ^{\prime
})-4\rbrack \lambda }{2(\eta -\eta ^{\prime })\hbar +6}+O(\lambda ^{3}),
\end{eqnarray*}
indicating that the Heisenberg algebras $\mathcal{H}_{\alpha }$, $\mathcal{H}%
_{\beta }$ are well-defined even at $\lambda =0$. The conjugates
of $\alpha (0)$ and $\beta (0)$ have to be introduced separately,
however, as follows. Let $Q_{\alpha }=\alpha (0),$ $Q_{\beta
}=\beta (0)$ and their conjugate operators $P_{\alpha }$,
$P_{\beta }$ are defined by the following relations:
\begin{eqnarray*}
\lbrack P_{\alpha },Q_{\alpha }\rbrack &=&1, \\
\lbrack P_{\beta },Q_{\beta }\rbrack &=&1, \\
\lbrack P_{\alpha },Q_{\beta }\rbrack &=&\lbrack P_{\beta },Q_{\alpha
}\rbrack =0.
\end{eqnarray*}
Now denoting
\begin{eqnarray*}
X_{a}(\lambda ) &=&\frac{1}{\hbar \lambda }\left( \mathrm{csch}\frac{\lambda
}{2\eta }\sinh \hbar \lambda +2\cosh \frac{\hbar \lambda }{2}+1\right) , \\
X_{b}(\lambda ) &=&\frac{1}{\hbar \lambda }\left( \mathrm{csch}\frac{\lambda
}{2\eta ^{\prime }}\sinh \hbar \lambda -2\cosh \frac{\hbar \lambda }{2}%
-1\right) , \\
Y_{a}(\lambda ) &=&Y_{b}(\lambda )=\frac{1}{\hbar \lambda }
\end{eqnarray*}
we can define
\begin{eqnarray*}
a(\lambda ) &=&X_{a}(\lambda )\alpha (\lambda )+Y_{a}(\lambda )\beta
(\lambda ), \\
b(\lambda ) &=&X_{b}(\lambda )\alpha (\lambda )+Y_{b}(\lambda )\beta
(\lambda ), \\
(\lambda &\neq &\mu )
\end{eqnarray*}
so that the corresponding commutation relations read:
\begin{eqnarray}
\lbrack a(\lambda ),a(\mu )\rbrack &=&-\frac{1}{\hbar ^{2}\lambda }\left( 1+%
\frac{\sinh \hbar \lambda }{\sinh \frac{\lambda }{2\eta }}\right) \delta
(\lambda +\mu ),  \label{coup1}\\
\lbrack b(\lambda ),b(\mu )\rbrack &=&-\frac{1}{\hbar ^{2}\lambda }\left( 1-%
\frac{\sinh \hbar \lambda }{\sinh \frac{\lambda }{2\eta ^{\prime }}}\right)
\delta (\lambda +\mu ), \label{coup2}\\
\lbrack a(\lambda ),b(\mu )\rbrack &=&\lbrack b(\lambda ),a(\mu )\rbrack =%
\frac{2}{\hbar ^{2}\lambda }\cosh \frac{\hbar \lambda }{2}\delta (\lambda
+\mu ),  \label{coup3}\\
(\lambda &\neq &\mu ) \nonumber.
\end{eqnarray}
These commutation relations are crucial for the construction of
the free boson representation for the algebra $\mathcal{A}_{\hbar
,\eta }(osp(2|2)^{(2)})$ and hence we give them a short name for
reference: $\mathcal{H}\lbrack a,b\rbrack $. Recall that we are
dealing with \emph{generic }deformation parameters, we do not
consider the specific values of the parameters at which the above
bosonic algebra becomes ill-defined (these include the points at
which $\frac{1}{2\eta }$ is a rational multiple of $\hbar $).

Before going into the details of the representation theory, we
have to specify
the structure of the Fock space on which the free bosonic algebra $\mathcal{H%
}\lbrack a,b\rbrack $ acts. Actually, there are infinite many ways
to realize the bosonic algebra $\mathcal{H}\lbrack a,b\rbrack $ in
terms of two commuting sets of Heisenberg algebras, and what we
outlined above is only one of infinitely many choices.

Denoting respectively by $\mathcal{F}_{\alpha }$ and $\mathcal{F}_{\beta }$
the Fock spaces for the Heisenberg algebras $\mathcal{H}_{\alpha }$ and $%
\mathcal{H}_{\beta }$, we see that the bosonic algebra
$\mathcal{H}\lbrack a,b\rbrack $ can be realized in a proper
subspace $\mathcal{F}\lbrack a,b\rbrack $ of $\mathcal{F}_{\alpha
}\otimes \mathcal{F}_{\beta }$ by actions of the form
\begin{eqnarray*}
\mathcal{F}_{\alpha }\otimes \mathcal{F}_{\beta } &\supset &\mathcal{F}%
\lbrack a,b\rbrack \ni
|v_{f_{1},...,f_{n},g_{1},...,g_{m}}\rangle=\int_{-\infty }^{-\epsilon
}d\lambda _{1}f(\lambda _{1})a(\lambda _{1})...\int_{-\infty }^{-\epsilon
}d\lambda _{n}f(\lambda _{n})a(\lambda _{n}) \\
&&\times \int_{-\infty }^{-\epsilon }d\mu _{1}g(\mu _{1})b(\mu
_{1})...\int_{-\infty }^{-\epsilon }d\mu _{m}g(\mu _{m})b(\mu _{m})|0\rangle
_{\alpha }\otimes |0\rangle _{\beta }, \\
&&\lambda _{1},...,\lambda _{n};\mu _{1},...,\mu _{m}<0,\ \forall n,m\in Z,
\\
&&0<\epsilon \rightarrow 0^{+}.
\end{eqnarray*}
Notice that the ordering of $a$ and $b$ in the above expression is
irrelevant, because all values of $\lambda _{1},...,\lambda _{n}$ and $\mu
_{1},...,\mu _{m}$ are negative. The Fock space thus described gives the
left action (or action onto the right) of the algebra $\mathcal{H}\lbrack
a,b\rbrack $. The Fock space which provides the right action (or action onto
the left) can be specified as the conjugation of the above,
{\em i.e.\ }$\mathcal{F}%
^{\ast }\lbrack a,b\rbrack $.

\vspace{1pt} It remains to specify the correlation functions for
operators acting on the Fock spaces $\mathcal{F}\lbrack a,b\rbrack $ and $%
\mathcal{F}^{\ast }\lbrack a,b\rbrack$, or using more precise mathematical
terminology, the pairing $\mathcal{F}\lbrack a,b\rbrack \otimes \mathcal{F}%
^{\ast }\lbrack a,b\rbrack \rightarrow C$. This is given by the
following three steps. First, we fix the normalization for the
vacuum vectors as follows,
\[
(_{\alpha }\langle 0|\otimes _{\beta }\langle 0|)(|0\rangle _{\alpha
}\otimes |0\rangle _{\beta })=1.
\]
Next, for any two vectors
\begin{eqnarray*}
\langle v_{f_{i}}| &=&_{\alpha }\langle 0|\otimes _{\beta }\langle
0|\int_{\epsilon }^{+\infty }d\lambda f_{i}(\lambda )X_{i}(\lambda ), \\
|v_{g_{j}}\rangle &=&\int_{-\infty }^{-\epsilon }d\mu g_{j}(\mu )X_{j}(\mu
)|0\rangle _{\alpha }\otimes |0\rangle _{\beta }
\end{eqnarray*}
where $X_{i,j}(\lambda )$ are operators acting on the Fock spaces $\mathcal{F%
}\lbrack a,b\rbrack $ and $\mathcal{F}^{\ast }\lbrack a,b\rbrack $
satisfying
\begin{eqnarray*}
X_{i,j}(\lambda )|0\rangle _{\alpha }\otimes |0\rangle _{\beta }
&=&0=_{\alpha }\langle 0|\otimes _{\beta }\langle 0|X_{i,j}(-\lambda
),\qquad (\lambda >0) \\
\lbrack X_{i}(\lambda ),X_{j}(\mu )\rbrack &=&x_{ij}(\lambda )\delta
(\lambda +\mu ),\qquad (x_{ij}(\lambda )\ \mathrm{regular\ at\ }\lambda =0),
\end{eqnarray*}
with $f_{i}(\lambda )$ and $g_{j}(\lambda )$ both analytic in a
small neighborhood of $\lambda =0$ except at $\lambda =0$ where
they have simple poles, we define the inner product as follows:
\[
\langle v_{f_{i}}|v_{g_{j}}\rangle =\int_{C}\frac{d\lambda \ln (-\lambda )}{%
2\pi i}f_{i}(\lambda )x_{ij}(\lambda )g_{j}(-\lambda ),
\]
where $C$ is an integration contour which goes from infinity to
zero above the positive real $\lambda $-axis, surrounding the
origin counterclockwise, and going to infinity again below the
positive real $\lambda $-axis. This particular kind of
regularization has already been used in \cite{KLP, f}. Last, for
`multi-particle' states like $\langle v_{f_{i_1},...,f_{i_k}}|$
and $|v_{g_{j_1},...,g_{j_k}}\rangle$ we apply the Wick theorem.

Having provided all the necessary tools for defining the free
boson representation, we now introduce the notation
\begin{eqnarray}
\varphi (u) &=&\int\limits_{\lambda \neq 0}d\lambda e^{i\lambda u}a(\lambda
), \\
\phi (u) &=&\int\limits_{\lambda \neq 0}d\lambda e^{i\lambda u}b(\lambda
),\label{phiphi}
\end{eqnarray}
where $\int\limits_{\lambda \neq 0}d\lambda $ means the
integration over the whole real $\lambda $ axis except the point
$\lambda =0$, {\em i.e.}
\[
\int\limits_{\lambda \neq 0}d\lambda =\lim_{\epsilon \rightarrow
0^{+}}\left( \int_{-\infty }^{-\epsilon }d\lambda +\int_{\epsilon }^{+\infty
}d\lambda \right) .
\]
We then have

\begin{proposition}
The following expressions give a free boson realization of the algebra $%
\mathcal{A}_{\hbar ,\eta }(osp(2|2)^{(2)})$ at $c=1$:
\begin{eqnarray*}
E(u) &=&e^{\gamma _{E}-\ln \eta }:\exp \left[ \frac{i\pi }{2}\left( \frac{1}{%
p_{\alpha }}P_{\alpha }+\frac{1}{p_{\beta }}P_{\beta }\right)
+\left( p_{\alpha }Q_{\alpha }+p_{\beta }Q_{\beta }\right) \right]
\exp \hbar \lbrack \varphi (u)\rbrack :\,, \\
F(u) &=&e^{\gamma _{E}-\ln \eta }:\exp \left[ \frac{i\pi }{2}\left( \frac{1}{%
p_{\alpha }}P_{\alpha }-\frac{1}{p_{\beta }}P_{\beta }\right)
+\left( p_{\alpha }Q_{\alpha }+p_{\beta }Q_{\beta }\right) \right]
\exp \hbar \lbrack \phi (u)\rbrack :\,, \\ H^{\pm }(u) &=&:\exp
\left[ i\pi \left( \frac{1}{p_{\alpha }}P_{\alpha }\right)
+2\left( p_{\alpha }Q_{\alpha }+p_{\beta }Q_{\beta }\right)
\right] \exp \hbar \lbrack \varphi (u\pm i\hbar /4)+\phi (u\mp
i\hbar /4)\rbrack :\,,
\end{eqnarray*}
where $p_{\alpha },$ $p_{\beta }$ are two arbitrary nonzero constants and $%
\gamma _{E}$ is the Euler constant $\gamma _{E}=$0.57721566....$\hfill
\square $
\end{proposition}

\vspace{1pt}The proof is by straightforward calculation using the
Fock space conventions above. The following formulae play crucial
roles:
\begin{eqnarray*}
\int_{C}\frac{d\lambda \ln (-\lambda )}{2\pi i\lambda }\frac{e^{-x\lambda }}{%
1-e^{-\lambda /\eta }} &=&\ln \Gamma (\eta x)+(\eta x-\frac{1}{2})(\gamma_E
-\ln \eta )-\frac{1}{2}\ln (2\pi ), \\
\Gamma \left( \frac{1}{2}-x\right) \Gamma \left( \frac{1}{2}+x\right) &=&%
\frac{\pi }{\cos \pi x}.
\end{eqnarray*}

\subsection{Free boson representations of $U_{q}(osp(2|2)^{(2)})$}

In this subsection we carry out the same practice as in the last subsection
without assuming the relation $1/\eta ^{\prime}=1/\eta +\hbar $. As we shall see,
this yields a representation of the algebra $%
U_{q}(osp(2|2)^{(2)})$ at $\gamma =q^{1/2}$. Below we give some of the
details.

\vspace{1pt}
We introduce
two Heisenberg algebras $\tilde{\mathcal{H}}_{\alpha }$, $\tilde{\mathcal{H}}%
_{\beta }$, not to be confused with those of the last subsection,
defined respectively by

\begin{eqnarray*}
\lbrack \alpha (\lambda ),\alpha (\mu )\rbrack &=&A(\lambda )\delta (\lambda
+\mu ), \\
\lbrack \beta (\lambda ),\beta (\mu )\rbrack &=&B(\lambda )\delta (\lambda
+\mu ), \\
\lbrack \alpha (\lambda ),\beta (\mu )\rbrack &=&0, \\
(\lambda &\neq &\mu )
\end{eqnarray*}
where $A(\lambda )$ and $B(\lambda )$ are given by

\begin{eqnarray*}
A(\lambda ) &=&\frac{\lambda }{4\cosh \frac{\hbar \lambda }{2}+2}, \\
B(\lambda ) &=&-\frac{\lambda \left( \mathrm{csch}^{2}\frac{\lambda }{2\eta }%
\sinh ^{2}\hbar \lambda +2\cosh \hbar \lambda +1\right) }{4\cosh \frac{\hbar
\lambda }{2}+2}.
\end{eqnarray*}
The zero mode operators remain as in the last subsection:
\begin{eqnarray*}
\lbrack P_{\alpha },Q_{\alpha }\rbrack &=&1, \\
\lbrack P_{\beta },Q_{\beta }\rbrack &=&1, \\
\lbrack P_{\alpha },Q_{\beta }\rbrack &=&\lbrack P_{\beta },Q_{\alpha
}\rbrack =0.
\end{eqnarray*}
Let
\begin{eqnarray*}
X_{a}(\lambda ) &=&\frac{1}{\lambda }\left( \mathrm{csch}\frac{\lambda }{%
2\eta }\sinh \hbar \lambda +2\cosh \frac{\hbar \lambda }{2}+1\right) , \\
X_{b}(\lambda ) &=&\frac{1}{\lambda }\left( \mathrm{csch}\frac{\lambda }{%
2\eta }\sinh \hbar \lambda -2\cosh \frac{\hbar \lambda }{2}-1\right) , \\
Y_{a}(\lambda ) &=&Y_{b}(\lambda )=\frac{1}{\lambda }
\end{eqnarray*}
and define
\begin{eqnarray*}
a(\lambda ) &=&X_{a}(\lambda )\alpha (\lambda )+Y_{a}(\lambda )\beta
(\lambda ), \\
b(\lambda ) &=&X_{b}(\lambda )\alpha (\lambda )+Y_{b}(\lambda )\beta
(\lambda ), \\
(\lambda &\neq &\mu )
\end{eqnarray*}
so that the corresponding commutation relations read:
\begin{eqnarray}
\lbrack a(\lambda ),a(\mu )\rbrack &=&-\frac{1}{\lambda }\left( 1+\frac{%
\sinh \hbar \lambda }{\sinh \frac{\lambda }{2\eta }}\right) \delta (\lambda
+\mu ), \label{comm1}\\
\lbrack b(\lambda ),b(\mu )\rbrack &=&-\frac{1}{\lambda }\left( 1-\frac{%
\sinh \hbar \lambda }{\sinh \frac{\lambda }{2\eta }}\right) \delta (\lambda
+\mu ), \\
\lbrack a(\lambda ),b(\mu )\rbrack &=&\lbrack b(\lambda ),a(\mu )\rbrack =%
\frac{2}{\lambda }\cosh \frac{\hbar \lambda }{2}\delta (\lambda +\mu ), \label{comm3} \\
(\lambda &\neq &\mu ). \nonumber
\end{eqnarray}
We introduce the free bosonic fields $\varphi (u)$ and $\phi (u)$
as in (\ref{phiphi}) but using the Heisenberg algebras described
in this subsection.

\begin{proposition} \label{prop4}
The following expressions give a free boson realization of the algebra $%
U_{q}(osp(2|2)^{(2)})$ currents (\ref{uq}-\ref{ant2}) at $\gamma =q^{1/2}$:
\begin{eqnarray*}
X^{+}(z) &=&e^{\gamma _{E}-\ln \eta }\frac{1}{\sqrt{2}z}:\exp \left[ \frac{%
i\pi }{2}\left( \frac{1}{p_{\alpha }}P_{\alpha }+\frac{1}{p_{\beta }}%
P_{\beta }\right) +\left( p_{\alpha }Q_{\alpha }+p_{\beta
}Q_{\beta }\right) \right]  \\ &\times &\exp \left[ \varphi \left(
\frac{\ln z}{2\pi \eta }\right) \right] :\,, \\
X^{-}(z) &=&e^{\gamma _{E}-\ln \eta }\frac{1}{\sqrt{2}z}:\exp \left[ \frac{%
i\pi }{2}\left( \frac{1}{p_{\alpha }}P_{\alpha }-\frac{1}{p_{\beta }}%
P_{\beta }\right) +\left( p_{\alpha }Q_{\alpha }+p_{\beta
}Q_{\beta }\right) \right]  \\ &\times &\exp \left[ \phi \left(
\frac{\ln z}{2\pi \eta }\right) \right] :\,,
\\
\psi ^{\pm }(z) &=&\frac{(2\pi )^{2}i\eta (q-q^{-1})}{\ln q}:\exp \left[
i\pi \left( \frac{1}{p_{\alpha }}P_{\alpha }\right) +2\left( p_{\alpha
}Q_{\alpha }+p_{\beta }Q_{\beta }\right) \right]  \\
&\times &\exp \left[ \varphi \left( \frac{\ln zq^{\pm 1/4}}{2\pi \eta }%
\right) +\phi \left( \frac{\ln zq^{\mp 1/4}}{2\pi \eta }\right)
\right] :\,,
\end{eqnarray*}
where $p_{\alpha },$ $p_{\beta }$ are two arbitrary nonzero constants, $%
\gamma _{E}$ is the Euler constant, and $q$ and $z$ are related to the
parameters $\eta $ and $\hbar $ via $q=e^{2\pi i\eta \hbar }$ and $z=e^{2\pi
\eta u}$.\hfill $\square $
\end{proposition}

\begin{remark}
We can set the parameter $\eta$ in the above free boson
representation to any fixed value without changing the
representation itself. In this sense, the parameter $\eta $ can be
thought of as redundant. Indeed, the algebra $U_q(osp(2|2)^{(2)})$
contains only one deformation parameter, and its representation
must necessarily contain no extra parameters.\hfill $\square
$
\end{remark}

\begin{remark}
According to the Footnote 1, the value of $\gamma =q^{1/2}$
corresponds to $c=1/2$ in \cite{YZ} . Thus the free boson
representation given in Proposition \ref{prop4} is somehow not the
most interesting one --- remember that for usual affine Lie
(super)algebras, only representations at integer values of $c$
have received attention, because only these representations are
known to be unitary.$\hfill \square $
\end{remark}

\begin{remark}
In \cite{YZ}, a free boson representation of the
$U_{q}(osp(2|2)^{(2)}) $ currents (\ref{uq}-\ref{ant2}) at $\gamma
=q$ ({\em i.e.\ }$c=1$) was given. However, that representation
does not have a well-defined limit as $q\rightarrow 1$ and hence
the authors of that paper called their representation a
`nonclassical' one. One can verify that our representation does
have a well-defined limit at this value of deformation
parameter.\hfill $\square $
\end{remark}

Although the representation of $U_{q}(osp(2|2)^{(2)})$ at $\gamma
=q^{1/2}$ is not of the most interesting class, we could, however,
use the same method to construct `interesting' representations.
Now, instead of (\ref{comm1}-\ref{comm3}), we introduce the
following set of bosonic algebras:
\begin{eqnarray*}
\lbrack a_{n},a_{m}\rbrack &=&-\frac{1}{n}(1+(-q)^{-n})\delta _{n+m,0}, \\
\lbrack b_{n},b_{m}\rbrack &=&-\frac{1}{n}(1-(-q)^{-n})\delta _{n+m,0}, \\
\lbrack a_{n},b_{m}\rbrack &=&\lbrack b_{n},a_{m}\rbrack =\frac{1}{n}%
(q^{n}+q^{-n})\delta _{n+m,0},
\end{eqnarray*}
together with the zero mode operators:
\begin{eqnarray*}
\lbrack P_{a},Q_{a}\rbrack &=&1, \\
\lbrack P_{b},Q_{b}\rbrack &=&1, \\
\lbrack P_{a},Q_{b}\rbrack &=&\lbrack P_{b},Q_{a}\rbrack =0.
\end{eqnarray*}
These bosonic commutation relations can also be realized in the
tensor product of Fock spaces of two commuting sets of Heisenberg
algebras, though we omit this here.

Defining
\begin{eqnarray*}
\varphi (z) &=&\sum_{n\neq 0}a_{n}z^{-n}+P_{a}\ln z+2Q_{a}+2Q_{b}, \\
\phi (z) &=&\sum_{n\neq 0}b_{n}z^{-n}-(P_{a}+P_{b})(\ln z+i\pi/2
)-2(Q_{a}-Q_{b}),
\end{eqnarray*}
we can easily prove

\begin{proposition}
The following bosonic expressions give a free boson representation of $%
U_{q}(osp(2|2)^{(2)})$ currents at $\gamma =q$:
\begin{eqnarray*}
E(z) &=&:\exp \varphi (z):, \\
F(z) &=&:\exp \phi (z):, \\
H^{\pm }(z) &=&:E(zq^{\pm 1/2})F(zq^{\mp 1/2}):.
\end{eqnarray*}
\hfill $\square $
\end{proposition}
Again, this free boson representation is well defined as $q\rightarrow 1$.

\section{Concluding Remarks}
The result in this article provides an example of a two-parameter
deformed quantum current algebra with the structure of an infinite
Hopf family of superalgebras and associated with a
non-simply-laced and twisted root system. To the authors'
knowledge, this is the first example of this kind and hence a
useful hint at the final classification of all such algebras.

As mentioned in the introduction, algebras with the structure of
an infinite Hopf family of (super)algebras have so far been
studied only in current realizations, and this limitation has made
understanding the relationship between these algebras and the
quasi-Hopf algebras a difficult task. Meanwhile, in most physical
applications, quantum-deformed algebras are best formulated in
terms of the Yang-Baxter relations (the celebrated `$RLL$'
relations). It is expected that if such a realization for the
algebras of the present kind can be achieved, then many of the
unsolved problems mentioned above could be treated more easily.
The $RLL$ realization of the two-parameter deformed algebras with
the structure of infinite Hopf families of algebras is currently
under investigation.

\vspace{1cm} \noindent \textbf{Acknowledgement}: L. Zhao would
like to thank the Department of Mathematics, University of York
for hospitality. Discussions with Xiang-Mao Ding and financial
support from the Royal Society of London, the UK PPARC and the
National Natural Science Foundation of China are also warmly
acknowledged.

\end{document}